# The Solutions of Nonlinear Heat Conduction Equation via Fibonacci&Lucas Approximation Method

*Zehra Pınar[a]   Turgut Öziş[b]*

[a]*Namık Kemal University, Faculty of Science & Letters, Department of Mathematics, , 59030 Tekirdağ, Turkey*

[b]*Ege University,* Science *Faculty, Department of Mathematics, 35100 Bornova-İzmir, Turkey*

Abstract: To obtain new types of exact travelling wave solutions to nonlinear partial differential equations, a number of approximate methods are known in the literature. In this study, we extend the class of auxiliary equations of Fibonnacci&Lucas type equations. The proposed Fibonnacci&Lucas approximation method produces many new solutions. Consequently, we introduce new exact travelling wave solutions of some physical systems in terms of these new solutions of the Fibonacci&Lucas type equation. In addition to using different ansatz, we use determine different balancing principle to obtain optimal solutions.

[b] Corresponding author: Zehra Pınar (zpinar@nku.edu.tr , +90282 2552706)

## 1. Introduction

The inspection of nonlinear wave phenomena physical systems is of great interest from both mathematical and physical points of view. In most cases, the theoretical modeling based on nonlinear partial differential equations (NLPDEs) can accurately describe the wave dynamics of many physical systems. The critical importance is to find closed form solutions for NLPDEs of physical significance. This could be a very complicated task and, in fact, is not always possible since in various realistic problems in physical systems. So, searching for some exact physically significant solutions is an important topic because of wide applications of NLPDEs in biology, chemistry, physics, fluid dynamics and other areas of engineering[2,3].

Since many of the most useful techniques in analysis are formal or heuristic the trend in recent years has also been to justify and provide the new procedures or methods rigorously[18]. Hence, over the past decades, a number of approximate methods for finding travelling wave solutions to nonlinear evolution equations have been proposed/or developed and furthermore modified [6-20]. The solutions to various evolution equations have been found by one or other of these methods. The technique of these methods consist of the solutions of the nonlinear evolution equations such that the target solutions of the nonlinear evolution equations can be expressed as a polynomial in an elementary function which satisfies a particular ordinary differential equation. Recently, to determine the solutions of nonlinear evolution equations, many exact solutions of various auxiliary equations have been utilized [22-27].

In this paper, we will examine the consequences of the choice of the Fibonacci&Lucas type equation for determining the solutions of the nonlinear evolution equation in consideration and more we search for additional forms of new exact solutions of nonlinear differential equations which satisfying Fibonacci&Lucas type equation(s). To obtain wave solutions of nonlinear partial differential equations via Fibonacci&Lucas transformation $n$ optimal index value is proposed.

## 2. Fibonacci & Lucas Polynomials

In this section, we determine Fibonacci&Lucas polynomials. Ordinary differential equations satisfied by two families of Fibonacci and Lucas polynomials are derived using identities which relate them to the generalized polynomials, and nonpolynomial solutions are deduced from corresponding solutions of the partial differential equations.

The single variable polynomials $F_n(1, z)$ and $L_n(1, z), n > 0$, with the properties

$$F_n(x, y) = x^n F_n(1, \eta),$$
$$L_n(x, y) = x^n L_n(1, \eta),$$

where $\eta = y/x^2$, are referred to as the Fibonacci and Lucas polynomials, respectively, by Doman and Williams [1]. Galvez and Devesa [4] have shown that they satisfy the ordinary differential equations

$$\eta(1+4\eta)\frac{d^2 F_n}{d\eta^2} - [(n-1) + 2(2n-5)\eta]\frac{dF_n}{d\eta} + (n-1)(n-2)F_n = 0, \tag{1}$$

$$\eta(1+4\eta)\frac{d^2 L_n}{d\eta^2} - [(n-1) + 2(2n-3)\eta]\frac{dL_n}{d\eta} + n(n-1)L_n = 0. \tag{2}$$

Using the earlier results, it can be shown that a second linearly independent solution of the Eq.(1) is $L_n(1,\eta)/\sqrt{|1+4\eta|}$, and a second linearly independent solution of the Eq.(2) is $\sqrt{|1+4\eta|}F_n(1,\eta)$.

The polynomials $F_n(\zeta,1)$ and $L_n(\zeta,1)$, also referred to as Fibonacci and Lucas polynomials by Hoggatt and Bicknell [5], are related to the generalized polynomials by

$$F_n(x, y) = y^{(n-1)/2} F_n(\zeta,1), L_n(x, y) = y^{n/2} L_n(\zeta,1),$$

where $\zeta = x/\sqrt{y}$, that they satisfy the ordinary differential equations

$$(4+\zeta^2)\frac{d^2 F_n}{d\zeta^2} + 3\zeta\frac{dF_n}{d\zeta} - (n^2 - 1)F_n = 0, \tag{3}$$

$$(4+\zeta^2)\frac{d^2 L_n}{d\zeta^2} + \zeta\frac{dL_n}{d\zeta} - n^2 L_n = 0, \tag{4}$$

which also have solutions $L_n(\zeta,1)/\sqrt{\zeta^2+4}$ and $\sqrt{\zeta^2+4}F_n(\zeta,1)$, respectively.

As seen above, Eqs(1)-(4) are depends on $n$ index value, so the solutions of Eqs(1)-(4) are different solutions for each $n$ index values. But, one of the important questions is 'Which $n$ index value produce an available basis for the given problem?' and the answer is given in the following section.

**3. Methodology**

The fundamental nature of the auxiliary equation technique is given by many authors in the literature [6-15] and it is applied to this new approximation . Let us have a nonlinear partial differential equation

$$P(u, u_x, u_t, u_{xx}, u_{xt}, u_{tt}, \cdots) = 0 \qquad (5)$$

and let by means of an appropriate transformation which is depended on Fibonacci&Lucas type equation(s), this equation is reduced to nonlinear ordinary differential equation

$$Q(u, u_\xi, u_{\xi\xi}, u_{\xi\xi\xi}, \cdots) = 0. \qquad (6)$$

For large class of the equations of the type (6) have exact solutions which can be constructed via finite series

$$u(\xi) = \sum_{i=0}^{N} a_i z^i(\xi) \qquad (7)$$

Here, $a_i, (i = 0,1, \ldots, N)$ are parameters to be further determined, $N$ is an integer fixed by a balancing principle and elementary function $z(\xi)$ is the solution of some ordinary differential equation referred to as the auxiliary equation[15,16,17,19, 21].

It is worth to point out that we happen to know the general solution(s), $z(\xi)$, of the auxiliary equation beforehand or we know at least exact analytical particular solutions of the auxiliary equation.

*The outline of the method*:

A) Define the solution of Eq.(6) by the ansatz in form of finite series in Eq.(7) where $a_i, (i = 0,1, \ldots, N)$ are parameters to be further determined, $N$ is an integer fixed by a balancing principle and elementary function $z(\xi)$ is the solutions of Eqs(1-4) be considered which are depended on $n$ index values which helps us to obtain wave solution for the given problem. But, for each $n$ index values, we obtain different solutions of auxiliary equations(1)-(4), so we need to determine $n$ index values to control the availability of the given problem. To determine $n$ index values of Fibonacci&Lucas function from the balancing principle we proposed a novel balancing given below

$$\begin{cases} \dfrac{\partial^m u}{\partial \zeta^m} \text{ is highest degree} \\ n \equiv \{-1, 0, 1\} \bmod(m) \end{cases} \qquad (8)$$

B) Substitute Eq.(7) into ordinary differential equation Eq.(6) to determine the parameters $a_i, (i = 0,1, \ldots, N)$ with the aid of symbolic computation.

C) Insert predetermined parameters $a_i$ and elementary function $z(\xi)$ of the auxiliary equation into Eq.(7) to obtain travelling wave solutions of the nonlinear evolution equation in consideration.

It is very apparent that determining the elementary function $z(\xi)$ via auxiliary equation is crucial and plays very important role finding new travelling wave solutions of nonlinear evolution equations. This fact, indeed, compel researchers to search for a novel auxiliary equations with exact solutions.

In this study, we use Fibonacci& Lucas differential equations are given above.

## 4. Travelling Wave Solutions of Nonlinear Heat Conduction Equation In Terms of Fibonacci&Lucas Equation Regarding As An Ansatz.

In this section, we consider the following nonlinear heat conduction equation

$$u_t - \left(u^2\right)_{xx} - pu + qu^2 = 0 \tag{9}$$

where $p, q$ are real constants. To use Fibonacci&Lucas approximation method, we consider the determined variables $\zeta, \eta$ instead of the wave variable. Although the wave transformation is not used to obtain the wave solution, it is obtained by Fibonacci&Lucas transformations which are determined by $\zeta, \eta$. The $\zeta, \eta$ variables carries Eq. (9) into the ordinary differential equation.

From the balancing principle, $N = 1$ is obtained. Therefore, the ansatz yields

$$U(\zeta) = g_0 + g_1 z(\zeta) \tag{10}$$

where $z(\zeta)$ may be determined by the solution of Eqs.(1)-(4).

**Case 1:** We consider $\eta = y/x^2$ transformation and Eq(1) is considered as an auxiliary equation,

The $\eta$ variables carries Eq. (9) into the ordinary differential equation

$$\frac{U'(\eta)}{x^2} - 8t^2 \frac{U'(\eta)^2}{x^6} - 8U(\eta)t^2 \frac{U''(\eta)}{x^6} - 12U(\eta)t \frac{U'(\eta)}{x^4} - pU(\eta) + qU(\eta)^3 = 0 \tag{11}$$

and here the ansatz is assumed as following

$$U(\eta) = g_0 + g_1 F_n(\eta) \tag{12}$$

Hence, substituting Eqs.(1) and(12) into Eq.(11) and letting each coefficient of $F_n(\eta)$ to be zero, we obtain algebraic equation system and solving the system by the aid of Maple 16, we can determine the coefficients:

i) $C_1 = C_1, C_2 = C_2, g_0 = -\dfrac{8\eta(2-3n+n^2)}{3qx^2(1+\eta)}$

$g_1 = \dfrac{64\eta^2 + 2144n\eta^3 + 256n^3\eta^3 - 1312n^2\eta^3 - 1088\eta^3 + 3qx^2 + 6qx^2\eta - 21qx^2\eta^2 - 48qx^2\eta^3 - 24qx^2\eta^4 + 64n^3\eta^2 - 160n^2\eta^2 + 32n\eta^2}{12qx^2(1+\eta)\eta(1-17\eta+2n+8n\eta)}$

ii) $C_1 = C_1, C_2 = \dfrac{px^2(1+\eta)^7}{5(-1-\eta+12g_1\eta-36g_1\eta^2+8\eta^2+8\eta^3)}$

$g_0 = -\dfrac{8\eta(2-3n+n^2)}{3qx^2(1+\eta)}$

Substituting the above coefficients into ansatz (12) with the solution of Lucas type equation, we obtain one of new solution of nonlinear heat conduction equation.

$$u(x,t) = -\dfrac{8\eta(2-3n+n^2)}{3qx^2(1+\eta)} + \dfrac{64\eta^2 + 2144n\eta^3 + 256n^3\eta^3 - 1312n^2\eta^3 - 1088\eta^3 + 3qx^2 + 6qx^2\eta - 21qx^2\eta^2 - 48qx^2\eta^3 - 24qx^2\eta^4 + 64n^3\eta^2 - 160n^2\eta^2 + 32n\eta^2}{12qx^2(1+\eta)\eta(1-17\eta+2n+8n\eta)}$$

$$\left( C_1 \, hypergeom\left(\left[\dfrac{9}{2}-2n-\dfrac{\sqrt{73-60n+12n^2}}{2}, \dfrac{9}{2}-2n+\dfrac{\sqrt{73-60n+12n^2}}{2}\right], [-3n+9], 1+\dfrac{t}{x^2}\right) \right.$$

$$\left. + C_2 \left(1+\dfrac{t}{x^2}\right)^{(-8+3n)} hypergeom\left(\left[-\dfrac{7}{2}+n-\dfrac{\sqrt{73-60n+12n^2}}{2}, -\dfrac{7}{2}+n+\dfrac{\sqrt{73-60n+12n^2}}{2}\right], [-7+3n], 1+\dfrac{t}{x^2}\right) \right)$$

$$u(x,t) = -\dfrac{8\eta(2-3n+n^2)}{3qx^2(1+\eta)} + g_1 \left( \dfrac{C_1 \, hypergeom\left(\left[\dfrac{9}{2}-2n-\dfrac{\sqrt{73-60n+12n^2}}{2}, \dfrac{9}{2}-2n+\dfrac{\sqrt{73-60n+12n^2}}{2}\right], [-3n+9], 1+\dfrac{t}{x^2}\right)}{+\dfrac{px^2(1+\eta)^7}{5(-1-\eta+12g_1\eta-36g_1\eta^2+8\eta^2+8\eta^3)}\left(1+\dfrac{t}{x^2}\right)^{(-8+3n)} hypergeom\left(\left[-\dfrac{7}{2}+n-\dfrac{\sqrt{73-60n+12n^2}}{2}, -\dfrac{7}{2}+n+\dfrac{\sqrt{73-60n+12n^2}}{2}\right], [-7+3n], 1+\dfrac{t}{x^2}\right)} \right)$$

To obtain suitable solution, $n=1$ is obtained from Eq(8). For the special values of parameters, the solutions are shown in Figure 1.

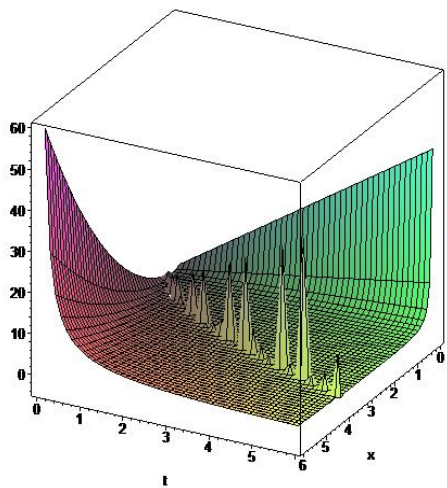 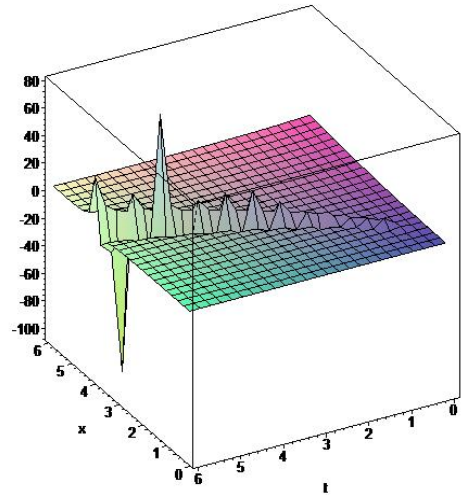

(a)                           (b)

**Figure 1.** (a) is for the first solution for $C_2 = 2, q = -1, C_1 = 0.5, n = 1$, (b) is for the second sloution where $p = 1, C_1 = 1, g_1 = -1, n = 1$,

If $n$ index values are chosen arbitrary, not using balancing formula Eq(8), the solutions of Eqs(1)-(4) don't construct a basis for the given equation. In that case, the behavior of solutions changes respect to $n$ index values. Therefore, to obtain travelling wave solutions for the given equation, we need to use the balancing formula Eq(8). For this example, if $n$ index values are taken as $n = 2, 3, 4, ...$, which are not satisfy the Eq(8), then for these index values we don't obtain solutions or these index values distort wave behavior of the obtained solutions. It is seen in the Figure 2.

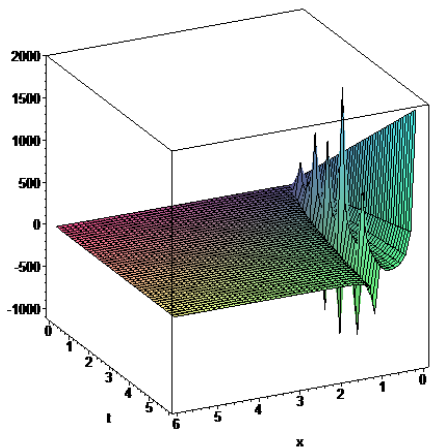 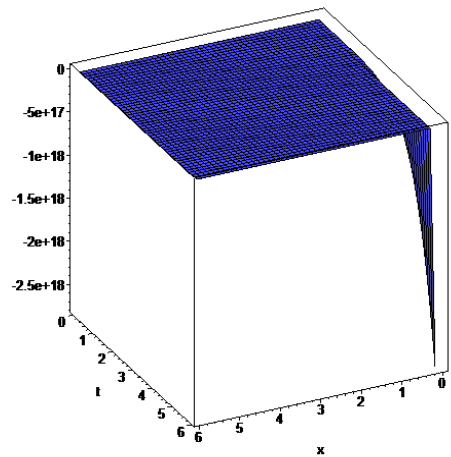

(a)                           (b)

**Figure 2.** (a) is for the $n=2$ index value, (b) is for the $n=4$ index value

**Case 2:** We consider $\zeta = x/\sqrt{y}$ transformation and Eq(4) is considered as an auxiliary equation. The $\zeta$ variables carries Eq. (9) into the ordinary differential equation

$$-\frac{U'(\zeta)x}{2t^{3/2}} - 2\frac{U'(\zeta)^2}{t} - 2U(\zeta)\frac{U''(\zeta)}{t} - pU(\zeta) - qU(\zeta)^3 = 0 \quad (13)$$

and here the ansatz is assumed as following from the balancing principle

$$U(\zeta) = g_0 + g_1 L_n(\zeta) \quad (14)$$

Hence, substituting Eqs.(4) and (14) into Eq.(13) and letting each coefficient of $L_n(\eta)$ to be zero, we obtain algebraic equation system and solving the system by the aid of Maple 16, we can determine the coefficients:

$$C_1 = C_1, C_2 = C_2, n = n, p = p, g_0 = 0, g_1 = g_1$$

and

$$C_1 = C_1, C_2 = \frac{\zeta\sqrt{-4-\zeta^2}\arctan\left(\frac{\zeta}{\sqrt{-4-\zeta^2}}\right)}{2\pi g_1}, n = n,$$

$$p = \frac{-\pi^2 - \arctan\left(\frac{\zeta\sqrt{-4-\zeta^2}}{4+\zeta^2}\right)^2 \zeta^2}{\arctan\left(\frac{\zeta\sqrt{-4-\zeta^2}}{4+\zeta^2}\right)^2 t(4+\zeta^2)}, g_0 = 0, g_1 = g_1$$

Substituting the above coefficients into ansatz (14) with the solution of Lucas type equation, we obtain new solutions of nonlinear heat conduction equation.

$$U(x,t) = x^n g_1 \left( C_1 \sin\left(\frac{1}{2}\sqrt{-8pt - 2px^2 - \frac{x^2}{t}}\arctan\left(\frac{x}{\sqrt{-4t-x^2}}\right)\right) + C_2 \cos\left(\frac{1}{2}\sqrt{-8pt - 2px^2 - \frac{x^2}{t}}\arctan\left(\frac{x}{\sqrt{-4t-x^2}}\right)\right) \right)$$

$$U(x,t) = x^n g_1 \left( C_1 \sin\left(n\arctan\left(\frac{x}{\sqrt{-4t-x^2}}\right)\right) + \frac{x\sqrt{-4-\frac{x^2}{t}}\arctan\left(\frac{x}{\sqrt{-4t-x^2}}\right)\cos\left(n\arctan\left(\frac{x}{\sqrt{-4t-x^2}}\right)\right)}{\pi g_1\sqrt{t}} \right)$$

To obtain optimal solution, $n=1$ is obtained from Eq(8). For the special values of parameters, the solutions are shown in Figure 3.

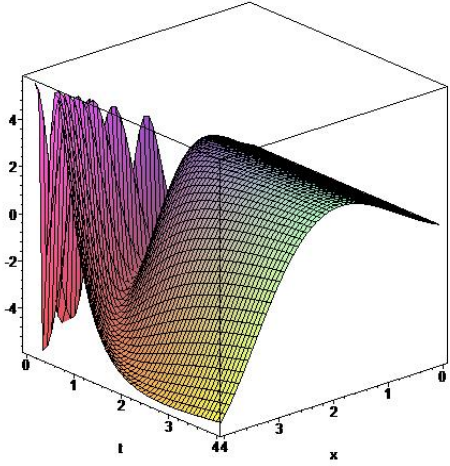
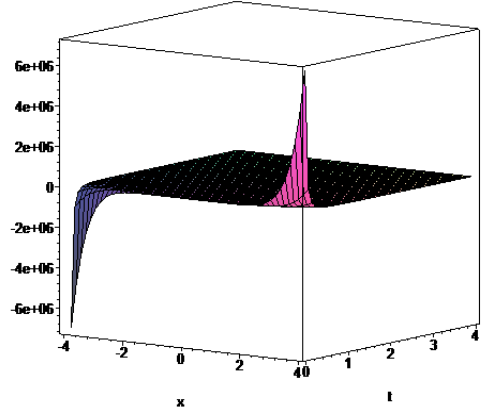

(a)                     (b)

**Figure 3:** (a) is the figure of the first solution, (b) is for the second solution for special values $C_1 = 1, p = 1, g_1 = 1, n = 1$.

**Case 3:** We consider $\eta = y/x^2$ transformation and Eq(2) is considered as an auxiliary equation,

The $\eta$ variables carries Eq. (9) into the ordinary differential equation

$$\frac{U'(\eta)}{x^2} - 8t^2 \frac{U'(\eta)^2}{x^6} - 8U(\eta)t^2 \frac{U''(\eta)}{x^6} - 12U(\eta)t\frac{U'(\eta)}{x^4} - pU(\eta) + qU(\eta)^3 = 0 \quad (11)$$

and here the ansatz is assumed as following

$$U(\eta) = g_0 + g_1 L_n(\eta) \quad (15)$$

Hence, substituting Eqs.(2) and (15) into Eq.(11) and letting each coefficient of $L_n(\eta)$ to be zero, we obtain algebraic equation system and solving the system by the aid of Maple 16, we can determine the coefficients:

i) $C_1 = 0, C_2 = -\dfrac{128\eta^2}{9qx^2(-1 - 3\eta + 5\eta^2 + 23\eta^3 + 24\eta^4 + 8\eta^5)}$

$g_0 = -\dfrac{8\eta n(n-1)}{3qx^2(1+\eta)}$

ii) $C_1 = 0,$

$C_2 = -\dfrac{4096\eta^3 - 144pqx^4\eta + 27pq^2x^6 g_1 + 81pq^2x^6 g_1\eta - 288pqx^4\eta^2 + 81pq^2x^6 g_1\eta^2 - 144pqx^4\eta^3 + 27pq^2x^6 g_1\eta^3}{18qx^2 g_1(1+\eta)^2(24qx^2\eta^4 - 448\eta^3 + 48qx^2\eta^3 + 84g_1qx^2\eta^3 + 144g_1qx^2\eta^2 + 21qx^2\eta^2 - 320\eta^2 - 6qx^2\eta + 60g_1qx^2\eta - 3qx^2)}$

$g_0 = -\dfrac{8\eta n(n-1)}{3qx^2(1+\eta)}$

iii) $p = \dfrac{2048\eta^4 n^3(2n^4 - 5n^3 - 1 + n - 35n\eta - 33n^3\eta + 8n^4\eta + 3n^2 + 51n^2\eta + 9\eta)}{18q^2 x^6 (1+\eta)^4 (-1+8\eta^2)}$

$g_1 = \dfrac{256 n^3 \eta^3 + 3qx^2 + 6qx^2\eta + 288n\eta^3 - 544n^2\eta^3 + 64n^3\eta^2 - 32n\eta^2 - 21qx^2\eta^2 - 48qx^2\eta^3 - 24qx^2\eta^4}{12qx^2\eta(1+\eta)(2n+8n\eta+1-9\eta)}$

$g_0 = -\dfrac{8\eta n(n-1)}{3qx^2(1+\eta)}$

Substituting the above coefficients into ansatz (15) with the solution of Lucas type equation, we obtain one of new solution of nonlinear heat conduction equation.

$u(x,t) = -\dfrac{8\eta n(n-1)}{3qx^2(1+\eta)} - \dfrac{128}{9} g_1 \eta^2 \left(1 + \dfrac{t}{x^2}\right)^{(-4+3n)} \dfrac{hypergeom\left(\left[-\dfrac{3}{2}+n-\dfrac{\sqrt{25-36n+12n^2}}{2}, -\dfrac{3}{2}+n+\dfrac{\sqrt{25-36n+12n^2}}{2}\right],[-3+3n],1+\dfrac{t}{x^2}\right)}{qx^2(-1-3\eta+5\eta^2+23\eta^3+24\eta^4+8\eta^5)}$

$u(x,t) = -\dfrac{8tn(n-1)}{3qx^4\left(1+\dfrac{t}{x^2}\right)} - \dfrac{4096\eta^3 - 144pqx^4\eta + 27pq^2x^6 g_1 + 81pq^2 x^6 g_1 \eta - 288pqx^4\eta^2 + 81pq^2x^6 g_1\eta^2 - 144pqx^4\eta^3 + 27pq^2 x^6 g_1\eta^3}{18qx^2 g_1(1+\eta)^2(24qx^2\eta^4 - 448\eta^3 + 48qx^2\eta^3 + 84g_1 qx^2\eta^3 + 144g_1 qx^2\eta^2 + 21qx^2\eta^2 - 320\eta^2 - 6qx^2\eta + 60g_1 qx^2\eta - 3qx^2)}\left(1+\dfrac{t}{x^2}\right)^{(-4+3n)}$

$hypergeom\left(\left[-\dfrac{3}{2}+n-\dfrac{\sqrt{25-36n+12n^2}}{2}, -\dfrac{3}{2}+n+\dfrac{\sqrt{25-36n+12n^2}}{2}\right],[-3+3n],1+\dfrac{t}{x^2}\right)$

$u(x,t) = -\dfrac{8tn(n-1)}{3qx^4\left(1+\dfrac{t}{x^2}\right)} + \dfrac{256n^3\eta^3 + 3qx^2 + 6qx^2\eta + 288n\eta^3 - 544n^2\eta^3 + 64n^3\eta^2 - 32n\eta^2 - 21qx^2\eta^2 - 48qx^2\eta^3 - 24qx^2\eta^4}{12qx^2\eta(1+\eta)(2n+8n\eta+1-9\eta)}$

$\left(\begin{array}{l} C_1 hypergeom\left(\left[\dfrac{5}{2}-2n-\dfrac{\sqrt{25-36n+12n^2}}{2}, \dfrac{5}{2}-2n+\dfrac{\sqrt{25-36n+12n^2}}{2}\right],[5-3n],1+\dfrac{t}{x^2}\right) \\ +C_2\left(1+\dfrac{t}{x^2}\right)^{(-4+3n)} hypergeom\left(\left[-\dfrac{3}{2}+n-\dfrac{\sqrt{25-36n+12n^2}}{2}, -\dfrac{3}{2}+n+\dfrac{\sqrt{25-36n+12n^2}}{2}\right],[-3+3n],1+\dfrac{t}{x^2}\right) \end{array}\right)$

To obtain optimal solution, $n=1$ is obtained from Eq(8). For the special values of parameters, the solutions are shown in Figure 4.

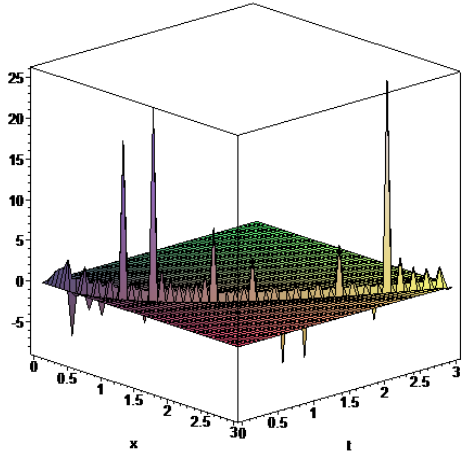
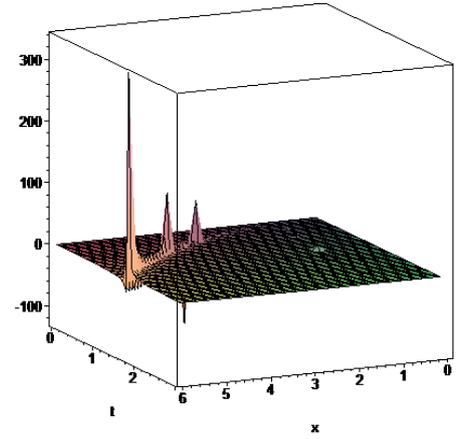

(a) (b)

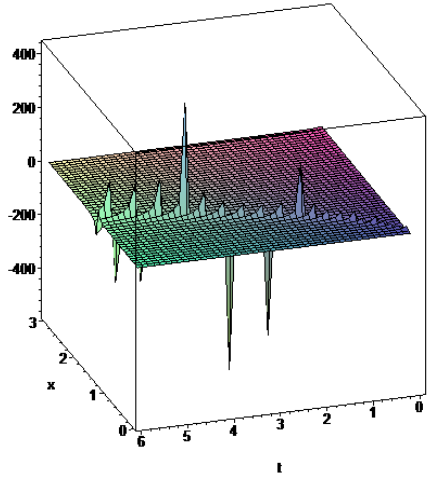

(c)

**Figure 4.** (a) is for the first solution for $C_2 = 1, q = 1, g_1 = 1, n = 1$, (b) is for the second sloution where $p = 1, q = 1, g_1 = 1, n = 1$, (c) is for the third solution $C_2 = 1, C_1 = 1, q = 1, n = 1$

**Case 4:** We consider $\zeta = x/\sqrt{y}$ transformation and Eq(3) is considered as an auxiliary equation, the $\zeta$ variables carries Eq. (9) into the ordinary differential equation

$$-\frac{U'(\zeta)x}{2t^{3/2}} - 2\frac{U'(\zeta)^2}{t} - 2U(\zeta)\frac{U''(\zeta)}{t} - pU(\zeta) - qU(\zeta)^3 = 0 \quad (13)$$

and here the ansatz is assumed as following

$$U(\eta) = g_0 + g_{-1}F_n^{-1}(\eta) \qquad (16)$$

Hence, substituting Eqs.(3) and (16) into Eq.(13) and letting each coefficient of $F_n(\eta)$ to be zero, we obtain algebraic equation system and solving the system by the aid of Maple 16, we can determine the coefficients:

$$C_1 = \frac{32\zeta^2 + 8\zeta\sqrt{4+\zeta^2} + 32\zeta^4 + 18\zeta^3\sqrt{4+\zeta^2} + 64\zeta^2 g_1 C_2 + 32\zeta g_1 C_2 \sqrt{4+\zeta^2} + 10\zeta^6 + 8\zeta^5\sqrt{4+\zeta^2} + 16\zeta^4 g_1 C_2 + 16\zeta^3 g_1 C_2 \sqrt{4+\zeta^2} + \zeta^8 + \zeta^7\sqrt{4+\zeta^2} + 32 g_1 C_2}{4g_1\left(4\zeta^2 + 2\zeta\sqrt{4+\zeta^2} + 4\zeta^4 + \zeta^3\sqrt{4+\zeta^2} + 2\right)}$$

Substituting the above coefficients into ansatz (16) with the solution of Lucas type equation, we obtain one of new solution of nonlinear heat conduction equation.

$$u(x,t) = g_1 \left( \frac{\left(\zeta + \sqrt{4+\zeta^2}\right)\left(32\zeta^2 + 8\zeta\sqrt{4+\zeta^2} + 32\zeta^4 + 18\zeta^3\sqrt{4+\zeta^2} + 64\zeta^2 g_1 C_2 + 32\zeta g_1 C_2\sqrt{4+\zeta^2} + 10\zeta^6 + 8\zeta^5\sqrt{4+\zeta^2} + 16\zeta^4 g_1 C_2 + 16\zeta^3 g_1 C_2\sqrt{4+\zeta^2} + \zeta^8 + \zeta^7\sqrt{4+\zeta^2} + 32g_1 C_2\right)}{\left(4g_1\left(4\zeta^2 + 2\zeta\sqrt{4+\zeta^2} + 4\zeta^4 + \zeta^3\sqrt{4+\zeta^2} + 2\right)\zeta + \dfrac{C_2}{\sqrt{4+\zeta^2}\left(\zeta + \sqrt{4+\zeta^2}\right)}\right)} \right)$$

To obtain optimal solution, $n = -1$ is obtained from Eq(8). For the special values of parameters, the solutions are shown in Figure 5.

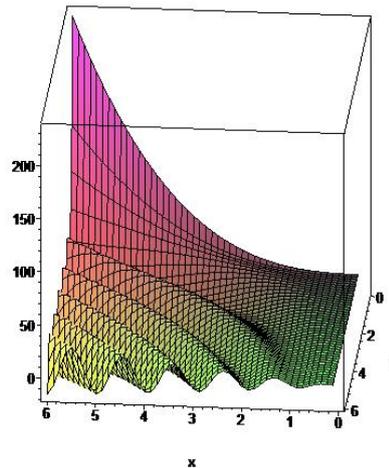

**Figure 5.** Graph of solution of Eq.(9) for $C_2 = 2, g_1 = \sin(xt), n = -1$

**Conclusion**

As it is seen, the key idea of obtaining new travelling wave solutions for the nonlinear equations is using the exact solutions of different types equations as an ansatz. By means of Fibonacci&Lucas type equation, the wave transformation is not used to obtain the wave solutions, the Fibonacci&Lucas transformation is used to obtain the wave solutions. But, the wave solutions are not obtained by the standard Fibonacci&Lucas transformations. For this reason, the $n$ optimal index value is proposed in Eq(8). Using the solutions of Fibonacci&Lucas type equation, we have successfully obtained a number of new exact periodic solutions of the nonlinear heat conduction equation by employing the solutions of the Fibonacci&Lucas type equation regarding as an auxiliary equation in proposed method.

In this letter, we have obtained new solutions of the nonlinear equation in hand using the Fibonacci&Lucas type equation (Eqs.(1-4)) for distinct cases. However, it is well known that the $n$ optimal index value produces new travelling wave solutions for many nonlinear problems. The presented method could lead to finding new exact travelling wave solutions for other nonlinear problems.

## 5. References

[1]     B. G. S. Doman & J. K. Williams. "Fibonacci and Lucas Polynomials.'1 *Math. Proc. Cambridge Philos. Soc.* 9.0 (1981):385-87.

[2]      L.Debnath, Nonlinear partial differential equations for scientists and Engineers (2$^{nd}$ ed.) Birkhauser, Boston (2005)

[3]     G.B.Whitham, A general approach to linear and nonlinear waves using a Lagrangian, *J. Fluid Mech.*,22(1965) 273-283

[4]      F. J. Galvez & J. S. Dehasa. "Novel Properties of Fibonacci and Lucas Polynomials." *MathProc. Cambridge Philos. Soc. 91* (1985): 159-6

[5]     V. E. Hoggatt, Jr., & M. Bicknell "Roots of Fibonacci Polynomials." *The Fibonacci Quarterly* 11.3 (1973):271-74

[6]     C. Yong, L. Biao, Z. Hong-Quing, Generalized Riccati equation expansion method and its application to Bogoyaylenskii's generalized breaking soliton equation, Chinese Physics. 12 (2003) 940–946.2004.

[7]     G. Cai, Q. Wang, J. Huang, A modified F-expansion method for solving breaking soliton equation, International Journal of Nonlinear science 2(2006) 122-128

[8]     X. Zeng, X. Yong, A new mapping method and its applications to nonlinear partial differential equations, Phy. Lett. A. 372 (2008) 6602–6607.

[9]     X. Yong, X. Zeng, Z. Zhang, Y. Chen, Symbolic computation of Jacobi elliptic function solutions to nonlinear differential-difference equations, Comput.Math. Appl. doi:10.1016/j.camwa2009.01.008.

[10]    W. X. Ma, T. Huang, Y. Zhang, A multiple Exp-function method for nonlinear differential equations and its applications, Phys. Scr. 82(2010) 065003(8pp)

[11]    T. Ozis, I. Aslan, Symbolic computation and construction of New exact traveling wawe solutions to Fitzhugh-Nagumo and Klein Gordon equation, Z. Naturforsch. 64a(2009) 15-20

[12]    T. Ozis, I. Aslan, Symbolic computations and exact and explicit solutions of some nonlinear evolution equations in mathematical physics, Commun. Theor. Phys., 51(2009) 577-580

[13]    Sirendaoreji, Auxiliary equation method and new solutions of Klein–Gordon equations,

Chaos, Solitons and Fractals 31 (2007) 943–950

[14]   B. Jang, New exact travelling wave solutions of nonlinear Klein–Gordon equations, Chaos, Solitons and Fractals 41 (2009) 646–654

[15]   X. Lv, S. Lai, Y.H. Wu, An auxiliary equation technique and exact solutions for a nonlinear Klein–Gordon equation, Chaos, Solitons and Fractals 41 (2009) 82–90

[16]   E. Yomba, A generalized auxiliary equation method and its application to nonlinear Klein–Gordon and generalized nonlinear Camassa–Holm equations, Physics Letters A 372 (2008) 1048–1060

[17]   J. Nickel, Elliptic solutions to a generalized BBM equation, Physics Letters A 364 (2007) 221–226

[18]   E.T. Whittaker, G.N. Watson, A Course of Modern Analysis, Cambridge Univ. Press, Cambridge, 1927.

[19]   E. Yomba, The extended Fan's sub-equation method and its application to KdV-MKdV, BKK and variant Boussinesq equations, Phys. Lett. A, 336(2005) 463-476

[20]   Z. Y. Yan, An improved algebra method and its applications in nonlinear wave equations,Chaos Solitons & Fractals 21(2004)1013-1021

[21]   M.A. Abdou, A generalized auxiliary equation method and its applications, Nonlinear Dyn 52 (2008) 95–102

[22]   Z. Pinar, T. Ozis, An Observation on the Periodic Solutions to Nonlinear Physical models by means of the auxiliary equation with a sixth-degree nonlinear term, Commun Nonlinear Sci Numer Simulat, 18(2013)2177-2187

[23]   N.K. Vitanov, Z. I. Dimitrova, K. N. Vitanov, On the class of nonlinear PDEs that can be treated by the modified method of simplest equation. Application to generalized Degasperis–Processi equation and b-equation, Commun Nonlinear Sci Numer Simulat 16 (2011) 3033–3044

[24]   RI. Ivanov, Water waves and integrability, Phil Trans R Soc A, 365(2007) 2267–80

[25]   L. Debnath, Nonlinear water waves, New York, Academic Press; 1994

[26]   RS. Johnson, The classical problem of water waves: a reservoir of integrable and nearly-integrable equations, J Nonlin Math Phys 10 (2003) 72–92

[27]   N. K. Vitanov, Modified method of simplest equation: Powerful tool for obtaining exact and approximate traveling-wave solutions of nonlinear PDEs, Commun Nonlinear Sci Numer Simulat 16 (2011) 1176–1185

**List of Figures**

**Figure 1.** (a) is for the first solution for $C_2 = 2, q = -1, C_1 = 0.5, n = 1$, (b) is for the second sloution where $p = 1, C_1 = 1, g_1 = -1, n = 1$,

**Figure 2.** (a) is for the $n = 2$ index value, (b) is for the $n = 4$ index value

**Figure 3:** (a) is the figure of the first solution , (b) is for the second solution for special values $C_1 = 1, p = 1, g_1 = 1, n = 1$.

**Figure 4.** (a) is for the first solution for $C_2 = 1, q = 1, g_1 = 1, n = 1$, (b) is for the second sloution where $p = 1, q = 1, g_1 = 1, n = 1$, (c) is for the third solution $C_2 = 1, C_1 = 1, q = 1, n = 1$

**Figure 5.** Graph of solution of Eq.(9) for $C_2 = 2, g_1 = \sin(xt), n = -1$